\documentclass[11pt]{amsart}

\usepackage{setspace}
\usepackage{amsfonts,amsopn,amssymb}
\usepackage{amsmath,latexsym,soul,cite,mathrsfs}
\usepackage{color,enumitem,graphicx}
\usepackage[colorlinks=true,urlcolor=blue,
citecolor=red,linkcolor=blue,linktocpage,pdfpagelabels,
bookmarksnumbered,bookmarksopen]{hyperref}
\usepackage[english]{babel}

\usepackage[left=2.7cm,right=2.7cm,top=3cm,bottom=3cm]{geometry}

\numberwithin{equation}{section}

\newtheorem{theorem}{Theorem}[section]
\newtheorem{lemma}[theorem]{Lemma}

\newtheorem{cor}[theorem]{Corollary}
\newtheorem{remark}[theorem]{Remark}

\newtheorem{claim}{Claim}[section]

\newtheorem{theoremletter}{Theorem}

\newcommand{\e}{\mathrm{e}}
\newcommand{\R}{{\mathbb R}}

\renewcommand{\epsilon}{\varepsilon}

\renewcommand{\theta}{{\vartheta}}

\renewcommand{\rightarrow}{\to}

\newcommand{\fim}{\hfill\rule{2mm}{2mm}}

\newcommand{\ud}{\mathrm{d}}
\newcommand{\N}{\mathbb{N}}
\newcommand{\dive}{\mathrm{div}}
\DeclareMathOperator{\supp}{\mathrm{supp}}
\DeclareMathOperator{\dist}{\mathrm{dist}}

\title[fractional  elliptic equation involving exponential growth ]{Concentration phenomena for  fractional  elliptic equations involving exponential critical growth \thanks{Research partially supported  INCTmat/CNPq/Brazil}}

\author[C. Alves]{Claudianor O. Alves}
\address[C. Alves]{Unidade Acad\^emica de Matem\'atica, 
Federal University of Campina  Grande,
\newline\indent
58429-900, Campina Grande-PB, Brazil}
\email{\href{mailto:coalves@dme.ufcg.edu}{coalves@dme.ufcg.edu}}

\author[J.M.\ do \'O]{Jo\~ao Marcos do \'O}
\address[J.M. do \'O]{Department of Mathematics,
Federal University of Para\'{\i}ba
\newline\indent
58051-900, Jo\~ao Pessoa-PB, Brazil}
\email{\href{mailto:jmbo@pq.cnpq.br}{jmbo@pq.cnpq.br}}

\author[O.H.\ Miyagaki]{Ol\'{i}mpio H.\ Miyagaki}
\address[O.\ Miyagaki]{Department of Mathematics, 
 Federal University of Juiz de Fora
\newline\indent 
36036-330  Juiz de Fora, Minas Gerais, Brazil}
\email{\href{mailto:olimpio@ufv.br}{olimpio@ufv.br}}

\thanks{Research partially supported by INCTmat/MCT/Brazil, CNPq, CAPES/Brazil.}
 

\date{\today}


\begin{document}

%%%%%%%%%%%%%%%%%%%%%%%%%%%%%%%%%%%%%%%%%%%%%%%%%%%%%%%%%%%%%%%%%%%%%%%%%%%%%%%%%%%%%%%%%%%%%%%%%%%%%%%%%%%%%%%%%%%%%%%%%%%%%%%%%%%%%%%%%%%%%%%%%%%%%%%%%%%%%%%%
%
%                               ABSTRACT
%
%%%%%%%%%%%%%%%%%%%%%%%%%%%%%%%%%%%%%%%%%%%%%%%%%%%%%%%%%%%%%%%%%%%%%%%%%%%%%%%%%%%%%%%%%%%%%%%%%%%%%%%%%%%%%%%%%%%%%%%%%%%%%%%%%%%%%%%%%%%%%%%%%%%%%%%%%%%%%%%%

%----------------------------------------------------------------------------------
%----------------------------------------------------------------------------------
\begin{abstract}
%----------------------------------------------------------------------------------
%----------------------------------------------------------------------------------
In this paper, we deal with the following singular perturbed fractional elliptic problem
$
\epsilon^{} (-\Delta)^{1/2}{u}+V(z)u=f(u)\,\,\, \mbox{in} \,\,\, \mathbb{R},
$
where $ (-\Delta)^{1/2}u$ is the square root of the Laplacian and $f(s)$ has exponential critical growth. Under suitable conditions on $f(s)$, we construct a localized bound state solution concentrating at an isolated component of the positive local minimum points of the potential of $V$ as $\epsilon$ goes to $0.$
%----------------------------------------------------------------------------------
%----------------------------------------------------------------------------------
\end{abstract}
%----------------------------------------------------------------------------------
%----------------------------------------------------------------------------------
%
\maketitle

%%%%%%%%%%%%%%%%%%%%%%%%%%%%%%%%%%%%%%

%\bigskip
%\begin{center}
%\begin{minipage}{8cm}
%\footnotesize
\tableofcontents
%\end{minipage}
%\end{center}

%%%%%%%%%%%%%%%%%%%%%%%%%%%%%%%%%%%%%%
%\bigskip

%----------------------------------------------------------------------------------
%----------------------------------------------------------------------------------
\section{{\Large Introduction}}
%----------------------------------------------------------------------------------
%----------------------------------------------------------------------------------
In this paper, we are concerned with existence and concentration of positive solutions for the following singular perturbed fractional elliptic problem
%%%%%%%%%%%%%%%%%%%%%%%%%%%%%%%%%%%%%%%%%
% The problem 
%%%%%%%%%%%%%%%%%%%%%%%%%%%%%%%%%%%%%%%%%
\begin{equation} \label{Pb}
\left\{
\begin{aligned} 
\varepsilon (- \Delta)^{1/2} u + V(z) u=f(u) & \quad \mbox{in} & \R, \\
u \in H^{1/2}(\R),\quad u > 0   & \quad \mbox{on} & \R,
\end{aligned} \tag{$P_{\epsilon}$}
\right. 
\end{equation}
\noindent where  $\epsilon $ is a small positive parameter, the potential $V$ is bounded away from zero, the nonlinearity $f(s)$ has exponential critical growth and $ (-\Delta)^{1/2}u$ is the square root of the Laplacian, which may be defined for smooth functions as
$$
{\mathcal F}((-\Delta)^{{1}/{2}}u)(\xi)=|\xi|{\mathcal F}(u)(\xi),
 $$ 
 where ${\mathcal F}$ is the Fourier transform, that is, 
 \[
 {\mathcal F}(\xi)=\frac{1}{\sqrt{2\pi}} \int_{\mathbb{R}} \e^{-i \xi \cdot x} \phi (x)  \, \ud x ,  
 \]
for functions $\phi$ in  the  Schwartz class.
Also, for sufficiently smooth $u$,  $(-\Delta)^{1/2}u$  can be equivalently represented , see \cite{nezza, Frank}, as
$$
(-\Delta)^{1/2} u = -\frac{1}{2\pi}\int_{\R}\frac{u(x+y)+u(x-y)-2 u(x)}{|y|^{2}} \ud y,
$$
and, by \cite[Propostion~3.6]{nezza},
\[
\|(-\Delta)^{1/4} u\|^2_{L^2}:=\frac{1}{2\pi}\int_{\R^2}\frac{(u(x)-u(y))^2}{|x-y|^{2}}\  \ud x \ud y, \quad \forall u\in  H^{1/2}(\R).
\]
Here $ H^{1/2}(\R)$ is the fractional Sobolev space
$$  H^{1/2}(\R)=\left\{ u \in L^2(\R): \ \|(-\Delta)^{1/4} u\|^2_{L^2}< \infty  \right\},    $$
endowed with the norm
$$\|u\|_{H^{1/2}}=\left(\|u\|^{2}_{L^2}+ \|(-\Delta)^{1/4} u\|^2_{L^2}\right)^{1/2}.$$

%-----------------------------------------------------------------------------------
%-----------------------------------------------------------------------------------
% Hypotheses on the potential $V$:
%-----------------------------------------------------------------------------------
%-----------------------------------------------------------------------------------
We suppose that the potential $V:\R \to \R$ is bounded and satisfies the following hypotheses:
\begin{description}
\item [($V_1$)] $ V$  is locally H\"older continuous and there  exists $V_0 >0$  such that
$$
V(z) \geq V_0, \quad \forall z \in \R, 
$$
\item [($V_2$)] there exists a bounded interval $\Lambda\subset \R$ such that
$$
V_0\equiv \inf_{\Lambda} V(z) < \min_{\partial \Lambda} V(z).
$$
\end{description}
%-----------------------------------------------------------------------------------
%-----------------------------------------------------------------------------------
% Hypotheses on the nonlinearity $f(s)$:
%-----------------------------------------------------------------------------------
%-----------------------------------------------------------------------------------

The function $f:\R \to \R$  satisfies the so-called Ambrosetti-Rabinowitz condition, introduced in \cite{ARMPT}, namely, 
\begin{equation}
\text{there exists $\theta >2$ with
           $0<\theta F(s)\leq sf(s)$ for all $s>0$, \,\,\, $F(s)=\int_{0}^{s}f(t)\; \ud t$.}
\label{AR}\tag{AR}
\end{equation}

In addition to the above condition we make the following assumptions on $f$:
\begin{description}
\item[(f1) ]  $f:\R\to\R^+$ is $C^1$ function with $f(s)=0$ if $s<0$.
\item[(f2) ]
$ \displaystyle
f(s)=o(s) \ \mbox{near origin}.
$
\item[(f3)]  $f(s)/{s}$ is an increasing function in $\R^+.$ 
\item[(f4)] There exist constants $p>2$ and $C_p>0$ such  that
$$
f(s)\geq C_p s^{p-1}\ \ \mbox{for all}\ \ s >0,
$$
where
$$
C_p>\left[\beta_p\left(\frac{2\theta}{\theta-2}\right)\frac{1}{\min\{1,V_0\}}\right]^{(p-2)/2},
$$
with
$$
\beta_p=\inf_{\mathcal{N}_0}\tilde{J}_0,
$$
$$
\mathcal{N}_0=\{v \in X^{1}(\mathbb{R}^{2}_+) \setminus \{0\}:\ \tilde{J}_{0}'(v)v=0\}
$$
and
$$
\tilde{J}_{0}(v)=\dfrac{1}{2}\int_{\mathbb{R}^{2}_+}|\nabla v|^2 \; \ud x \ud y +\dfrac{1}{2}\int_{\mathbb{R}}V_0|v(x,0)|^2\; \ud x-\dfrac{1}{p}\int_{\mathbb{R}}|v(x,0)|^p \; \ud x .
$$
where $X^{1}(\mathbb{R}^{2}_+)$  is defined in \eqref{X1}.
\end{description}

We are interested in bound state solution of \eqref{Pb} (solution with finite energy),
when $f$ has the maximal growth which allow us to treat problem \eqref{Pb} variationally in the fractional Sobolev space $H^{1/2}(\R)$ motivated by the following Trudinger-Moser
type inequality due to T.~Ozawa \cite{Ozawa}.
%-----------------------------------------------------------------------------------
%-----------------------------------------------------------------------------------
%	{Ozawa inequality}
%-----------------------------------------------------------------------------------
%-----------------------------------------------------------------------------------
\begin{theoremletter}
\label{moser}
 There exists $0 < \omega \leq \pi$ such that, for all $\alpha \in (0, \omega)$, there exists $H_{\alpha}>0$ with
\begin{equation}\label{I1}
\int_{\R} (\e^{\alpha u^2}-1) \, \ud x \leq H_{\alpha}\|u\|^{2}_{L^2},
\end{equation}
for all $u\in H^{1/2}(\R)$ with $\|(-\Delta)^{1/4} u\|^2_{L^2}\leq 1$.
\end{theoremletter}
\noindent
%-----------------------------------------------------------------------------------
%-----------------------------------------------------------------------------------
%	{Critical growth}
%-----------------------------------------------------------------------------------
%-----------------------------------------------------------------------------------
In view of \eqref{I1}, we say that $f$ has exponential critical growth at $+\infty,$ if there exist $\omega \in (0,\pi)$ and 
 $\alpha_0 \in (0, \omega),$
such that
$$
\lim_{s\to +\infty}\frac{f(s)}{\e^{\alpha s^2}}=0, \ \forall \alpha > \alpha_0,\quad \mbox{and} \quad  
\lim_{s\to +\infty}\frac{f(s)}{\e^{\alpha s^2}}=+\infty, \ \forall \alpha < \alpha_0.
$$
%-----------------------------------------------------------------------------------
%-----------------------------------------------------------------------------------
\subsection{Statement of the main the result}
%-----------------------------------------------------------------------------------
%-----------------------------------------------------------------------------------

The following theorem contains our main result:

\begin{theorem}
\label{resultado} Assume $(V_1),(V_2),$ $(AR)$, and $(f1)-(f4)$ hold. Then there exists $\epsilon_0 >0$  such that for $\epsilon \in (0,\epsilon_0)$, problem  \eqref{Pb} possesses a positive bound state solution $u_{\epsilon}(z)$ verifying the following conditions
\begin{description}
\item[I)] $u_\epsilon$ has at most one local (hence global) maximum $z_\epsilon $ in $ \R$ and $z_\epsilon \in I,$
\item[II)]$\displaystyle \lim_{\epsilon \to 0^+} V(z_\epsilon)=V_1=\inf_{I}V, $
\end{description}
\end{theorem}

This result extends to the nonlocal case the main result  in \cite{OS}. The proof is made combining Ozawa inequality \cite{Ozawa} with Del Pino and Felmer \cite{Pino} truncation argument and a recent approach developed in Alves and Miyagaki \cite{penalizacao}. In \cite{Secchi,Shang, CW} were established  existence results in nonlocal situation, while in \cite{Moustapha,ShangZhang, Chen,Davila} a concentration phenomena were proved imposing a global condition in $V.$

\begin{remark}
\begin{itemize}
\item[(1)] We recall that the condition (\ref{AR}) impose some superquadratic growth condition on the nonlinearity $F$.
\item[(2)] The condition $(f4)$ appeared first in \cite{Cao},  then for instance in \cite{AdoOM} and  \cite{OS}.  For the non-local situation it was used, e.g., in\cite{DMS}.
\item[(3)]  Critical growth of Trudinger-Moser type was  used in \cite{DMR}, also in \cite{OS,AlvesNovo,AdoOM}. In \cite{Antonio}  and in \cite{DMS} were used the Ozawa inequality to discuss nonlocal problem in bounded and unbounded domain, respectively.
\item[(4)] Notice that, if $f(s)$ has exponential critical growth, instead of assumption $(f4)$, it is enough to assume that there exist $p>2$ and $\mu >0$ such that
\[
\liminf_{s\rightarrow 0^+} \frac{f(s)}{s^{p-1}}  \geq \mu .
\]
\end{itemize}
\end{remark}

Throughout the paper, unless explicitly stated, the symbol $C$ will always denote
a generic positive constant, which may vary from line to line. 
%-----------------------------------------------------------------------------------
%-----------------------------------------------------------------------------------
\subsection{Outline} The sequel of the paper is organized as follows. The next section contains some technical results, which are crucial tools to prove our main theorem. In Sect.~\ref{Caffarelli}, we adapt a method due to L.~Caffarelli and L. Silvestre to obtain a local realization of the fractional Laplacian via a Dirichlet-to-Neummann operator. As a consequence of this argument we transform our nonlocal Problem~\eqref{Pb} into one local problem defined on the upper half plane \eqref{NPS}. Using variational techniques combined with Del~Pino and Felmer truncation argument we give the proof of Theorem~\ref{resultado} in Sect.~\ref{proof-Thm}, 
%-----------------------------------------------------------------------------------
%-----------------------------------------------------------------------------------

%-----------------------------------------------------------------------------------
%-----------------------------------------------------------------------------------
\section{Preliminary results}
%-----------------------------------------------------------------------------------
%-----------------------------------------------------------------------------------
In this section we collect preliminary facts for future reference. First
of all, let us set the standard notations to be used in the paper. We
denote the upper half-space in $\mathbb{R}^2$ by
$
\mathbb{R}^2_{+}=\{(x,y)\in \mathbb{R}^2 : y>0\}
$. In the sequel, $X^{1}(\R^{2}_{+})$ denotes the completion of $C^{\infty}_{0}(\overline{\R^{2}_{+}})$  with relation  to the norm $\|v\|_\epsilon $,
\begin{equation}\label{X1}
\begin{aligned}
X^{1}(\R^{2}_{+}) := & \overline{C^{\infty}_{0}(\overline{\R^{2}_{+}})}^{\|\cdot \|_\epsilon}, \quad \mbox{where} \\
\|v\|_\epsilon := & \Big(\int_{ \R^{2}_{+}}|\nabla v(x,y)|^2 \,\ud x\ud y + \int_{\R}V(\epsilon x)|v(x,0)|^2 \,\ud x\Big)^{1/2}.
\end{aligned}
\end{equation}
Moreover, we denote by $\|\,\,\,\,\|$ the usual norm in $X^{1}(\R^{2}_{+})$, that is, 
$$
\|v\|= \Big(\int_{ \R^{2}_{+}}|\nabla v(x,y)|^2 \,\ud x\ud y + \int_{\R}|v(x,0)|^2 \,\ud x\Big)^{1/2}.
$$
Since the potential $V$ is bounded from above and below, it is easy to see that $\|\,\,\,\|_\epsilon$ and $\|\,\,\,\|$ are equivalent norms in  $X^{1}(\R^{2}_{+})$ with
\begin{equation} \label{norma}
\min\{1,V_0\}\|v\| \leq \|v\|_\epsilon \leq \min\{1,|V|_{\infty}\}\|v\|, \quad \forall v \in  X^{1}(\R^{2}_{+}). 
\end{equation}

Using the above definition, we see  that if $v \in X^{1}(\mathbb{R}^{2}_+)$, then $u(x)=v(x,0)$ belongs to $H^{1/2}(\mathbb{R})$ and 
$$
\|v\|=\|u\|_{H^{1/2}}.
$$
Since $H^{1/2}(\mathbb{R})$ is continuously embedded into $L^{q}(\mathbb{R})$ for all $ q \geq 2$, c.f. \cite[Theorem~6.9]{nezza}, it follows that
$X^{1}(\mathbb{R}^{2}_+)$ is also continuously embedded into $L^{q}(\mathbb{R})$ for all $ q \geq 2$.  Moreover, the embedded 
$$
X^{1}(\mathbb{R}^{2}_+) \hookrightarrow L^{q}(A)
$$
are compact for any bounded mensurable set $A \subset \mathbb{R}$.  See \cite[Proposition~3.6]{Frank} also \cite[Remark~2.1]{DMS}.

Our first lemma is an important Trudinger-Moser  inequality on $X^{1}(\mathbb{R}^{2}_+)$, which was proved in \cite[Lemma 2.4]{DMS}.

\begin{lemma}
\label{boundmoser}
Let $(v_n)\subset X^{1}(\mathbb{R}^{2}_+)$ be a bounded sequence and assume $\displaystyle \sup_{n\in\N}\|v_n\|^{2}=M$. Then
$$
\sup_{n\in\N}\int_{\R} (e^{\alpha |v_n(x,0)|^2}-1) \,\ud x<\infty,\quad\text{for every $0<\alpha<\frac{\omega}{M^2}$};
$$
In particular, if $M\in (0,1)$, there exists $\alpha_M<\omega$ such that
$$
\sup_{n\in\N}\int_{\R} (e^{\alpha_M |v_n(x,0)|^2}-1) \,\ud x<\infty.
$$
\end{lemma}

Using the above lemma, we are able to prove some technical lemmas. The first of them is crucial in the study the behavior of Palais-Smale sequences. 
\begin{lemma} \label{alphat11} Let $(v_{n})$ be a sequence in $X^{1}(\R^{2}_{+})$ with 
	\begin{equation}\label{forro}
	\limsup_{n \to +\infty} \|v_n \|^{2}  < 1.
	\end{equation}
	Then, there exist  $t > 1$ sufficiently close to $ 1$ and $C > 0$  satisfying
	\[
	\int_{\mathbb{R}}\left(e^{\omega|v_n(x,0)|^{2}} - 1 \right)^t \; \ud x \leq C, \,\,\,\,\forall\, n \in \mathbb{N}. 
	\]

\end{lemma}
\noindent {\bf Proof.} \, 
Using \eqref{forro}
there are $m >0$ and $n_0\in \mathbb{N}$ verifying 
\[
||v_n||^{2} < m < 1,
\,\,\,\,\forall\,  n \geq n_0.
\]
Fix  $t > 1$ sufficiently close to $ 1$ and $ \beta> t$ satisfying $\beta m < 1. $ Then, there exists $C=C(\beta)>0$ such that
$$
\int_{\mathbb{R}}\left(e^{\omega|v_n(x,0)|^{2}} -1 \right)^t \; \ud x   
\leq C\int_{\mathbb{R}}\left(e^{\beta m  \omega(\frac{|v_n(x,0)|}{||v_n||})^{2}} - 1 \right) \; \ud x,
$$
for every  $n \geq n_0$.  Hence, by Lemma \ref{boundmoser},
$$
\int_{\mathbb{R}}\left(e^{\omega |v_n(x,0)|^{2}} -1 \right)^t \; \ud x \leq C_1 \, \,\,\,\,\forall\,  n \geq n_0,
$$
for some positive constant $C_1$. Now, the lemma follows fixing
$$
C=\max\left\{C_1, \int_{\mathbb{R}}\left(e^{\omega |v_{1}|^{2}} -1 \right)^t \; \ud x,....,\int_{\mathbb{R}}\left(e^{\omega |v_{n_0}|^{2}} -1 \right)^t \; \ud x \right\}.
$$
\fim

\begin{cor} \label{Convergencia em limitados} Let $(v_{n})$ be a sequence in $X^{1}(\mathbb{R}^{2}_+)$ satisfying \eqref{forro}.
If $v_n \rightharpoonup v$ weakly in $X^{1}(\mathbb{R}^{2}_+)$ and $v_n(x,0) \to v(x,0)$ a.e in $\mathbb{R}$, as $n \rightarrow \infty, $  then, 
\begin{equation} \label{EQUAT1}
F(v_n(x,0)) \to F(v(x,0)) \quad \mbox{in} \quad L^{1}(-R,R), 
\end{equation}
\begin{equation} \label{EQUAT1'}
f(v_n(x,0))v_n(x,0) \to f(v(x,0))v(x,0) \quad \mbox{in} \quad L^{1}(-R,R) 
\end{equation}
and
\begin{equation} \label{EQUAT2}
\int_{-R}^{R}f(v_n(x,0))\phi(x,0) \to \int_{-R}^{R}f(v(x,0))\phi(x,0), 
\end{equation}
as $n \rightarrow \infty, $ for all $\phi \in X^{1}(\mathbb{R}^{2}_+)$ and $R>0$. 
	
\end{cor}
\noindent {\bf Proof.} By $(f_1)$, for each $\beta >1$ and $\alpha > \alpha_0,$ there is $C>0$ such that
$$
|F(s)| \leq C(|s|^{2}+ (e^{\alpha \beta |s|^{2}}-1)) \quad \forall s \in \mathbb{R},
$$
from where it follows that, 
\begin{equation} \label{Domina}
|F(v_n(x,0))| \leq C(|v_n(x,0)|^{2}+(e^{\alpha \beta|v_n(x,0)|^{2}}-1)), \quad \forall n \in \mathbb{N}.
\end{equation}
Setting
$$
h_n(x)=C(e^{\alpha_0 \beta  |v_n(x,0)|^{2}}-1),
$$
we can fix $\beta, q>1$ sufficiently close to $1$ and $\alpha $ sufficiently close to $ \alpha_0$ such that 
$$
h_n \in L^{q}(\mathbb{R}) \quad \mbox{and} \quad \sup_{n \in \mathbb{N}}|h_n|_q<+\infty, 
$$
which is an immediate consequence of Lemma \ref{alphat11}. Therefore, up to subsequence, we derive that
$$
h_n \rightharpoonup h=C(e^{\alpha_0 \beta |v(x,0)|^{2}}-1) \quad \mbox{ weakly in} \quad L^{q}(\mathbb{R}), \ \mbox{as}\ n \to \infty.
$$
Since  $h_n, h \geq 0$, the last limit yields 
$$
h_n \to h \quad \mbox{in} \quad L^{1}(-R,R), \quad \forall R>0, , \ \mbox{as}\ n \to \infty.
$$
On the other hand, we know that
$$
v_n(\cdot,0) \to v(\cdot,0) \quad \mbox{in} \quad L^{2}(-R,R), \ \mbox{as}\ n \to \infty.  
$$
Gathering the above limits  with (\ref{Domina}), we get 
$$
F(v_n(x,0)) \to F(v(x,0)) \quad \mbox{in} \quad L^{1}(-R,R), \quad \forall R>0, \ \mbox{as}\ n \to \infty.
$$
The limits (\ref{EQUAT1'}) and (\ref{EQUAT2}) follow with the same type of arguments. 
\fim

The next lemma is a Lions type result, which is crucial in our approach. Since it follows with the same arguments found in  C.~Alves, J.~M.~do~\'O and O.~Miyagaki  \cite[Proposition 2.3]{AdoOM}, we will omit its proof.

\begin{lemma} \label{lions} Let $(v_n) \subset X^{1}(\mathbb{R}^{2}_+)$ be a sequence  with 
$$
\limsup_{n \to +\infty} \|v_n \|^{2}  < 1.
$$
If there is $R>0$ such that
$$
\lim_{n \to +\infty}\sup_{z \in \mathbb{R}}\int_{z-R}^{z+R}|v_n(x,0)|^{2}\, \ud x=0,
$$
then
$$
\lim_{n \to +\infty}\int_{\mathbb{R}}F(v_n(x,0))\,\ud x=\lim_{n \to +\infty}\int_{\mathbb{R}}f(v_n(x,0))v_n(x,0)\,\ud x=0.
$$
\end{lemma}

%-----------------------------------------------------------------------------------
%-----------------------------------------------------------------------------------
\section{The Caffarelli and Silvestre's method}\label{Caffarelli}
%-----------------------------------------------------------------------------------
%-----------------------------------------------------------------------------------

First of all, using the change variable $u(x)=v(\epsilon x)$, it is possible to prove that Problem \eqref{Pb}  is equivalent to the problem

\begin{equation} \label{Pa}
\left\{
\begin{aligned} 
 (- \Delta)^{1/2} u + V( \varepsilon z) u=f(u) & \quad \mbox{in} & \R, \\
u \in H^{1/2}(\R),\quad u > 0   & \quad \mbox{on} & \R.
\end{aligned} \tag{$P'_{\epsilon}$}
\right. 
\end{equation}
Hereafter, to get a solution for $(P'_{\epsilon})$, we will use a method due to L.~Caffarelli and L.~Silvestre in \cite{caffarelli}, more exactly, due to R.~Frank and E.~Lenzmann \cite{Frank} for whole line. In the seminal above  papers, were developed a local interpretation of the fractional Laplacian given in $\R$ by considering a Dirichlet to Neumann type operator in the domain
$\R^{2}_{+}=\{(x, t) \in \R^{2} : t > 0\}.$ A similar extension, in a bounded domain, see for instance, \cite{barrios, colorado, cabre}.
 For $u \in H^{1/2}(\R),$ the solution $w \in 
X^{1}(\R^{2}_{+})$ of  
 \begin{equation}
 \left\{ \begin{array}{rcl}
 -\dive ( \nabla w)=0 & \mbox{in}&  \R^{2}_{+}\noindent\\
  w=u& \mbox{on}& \R \times\{0\}\noindent
 \end{array}\right. 
 \end{equation}
is called $1/2$-harmonic  extension $w=E_{1/2}(u)$ of $u$ and it is proved in \cite{caffarelli} that
$$
\lim_{y \to 0^+} \frac{\partial w}{\partial y}(x,y)=-(-\Delta)^{1/2}u(x).
$$
%where  
%$$
%k_{s}=2^{1-2s}\Gamma(1-s)/\Gamma (s).
%$$ 

To get a solution for the nonlocal Problem \eqref{Pa}, we  will study the existence of solution for the local problem defined on the upper half plane
\begin{equation}\label{NPS}
\left\{ 
\begin{aligned}
-\dive ( \nabla w)= &\, 0 & \quad \mbox{in} & \quad  \R^{2}_{+}\noindent\\
 -\dfrac{\partial w}{\partial \nu}= & -V(\epsilon x)w + f(w)  &\quad \mbox{on}& \quad \R \times\{0\},\noindent
\end{aligned} \tag{$LP_{\epsilon}$}
\right. 
\end{equation}
where 
$$
\frac{\partial w}{\partial \nu}=\lim_{y \to 0^+} \frac{\partial w}{\partial y}(x,y),
$$ 
since if $w$ is a solution for the above problem, the function $u(x)=w(x,0)$ is a solution for \eqref{Pa}. 

Associated with \eqref{NPS}, we have the $J_\epsilon:X^{1}(\mathbb{R}^{2}_+) \to \mathbb{R}$ defined by
\begin{equation} \label{E1}
J_\epsilon (v)=\frac{1}{2}\int_{\mathbb{R}_{+}^{2}}|\nabla v|^{2}\,\ud x \ud y+
\frac{1}{2}\int_{\mathbb{R}}V(\epsilon x)|v(x,0)|^{2}\,\ud x-\int_{\mathbb{R}}F(v(x,0))\,\ud x,
\end{equation}
which is $C^{1}(X^{1}(\mathbb{R}^{2}_+),\mathbb{R})$ with derivative given by
\begin{eqnarray} \label{E2}
J'_\epsilon (v)\phi &=& \frac{1}{2}\int_{\mathbb{R}_{+}^{2}}\nabla v . \nabla \phi\,\ud x\ud y \\ &&+
\frac{1}{2}\int_{\mathbb{R}}V(\epsilon x) v(x,0) \phi(x,0) \,\ud x-\int_{\mathbb{R}}f(v(x,0))\phi (x,0)\,\ud x, \;\; \forall \phi \in X^{1}(\mathbb{R}^{2}_+).  \nonumber
\end{eqnarray}

We would like point out that $u$ is a solution of $(P'_{\epsilon})$ if, and only if, $u=v(x,0)$ for all $x \in \mathbb{R}$, for some critical point $v$ of $J_\epsilon$.

In what follows,  we will not work directly with functional $J_\epsilon$, because we have some difficulties to prove that it verifies the $(PS)$ condition. Hereafter, we will use the same approach explored in \cite{Pino}, modifying the nonlinearity of a suitable way. 
The idea is the following:  

First of all, without loss of generality, we will assume that
\begin{equation} \label{V_00}
0 \in \Lambda \quad \mbox{and} \quad V(0)=V_0=\inf_{x \in \R}V(x).
\end{equation}
We recall that in assumption $ (f1) $ we imposed that $f(t)=0,\; \forall t \leq 0 $,  because we are looking for positive solutions. 
Moreover, let us choose $k > { 2 \theta}/{(\theta -2)}$ and $a>0$ verifying
$$
\frac{f(a)}{a}=\frac{V_0}{k},
$$
where $V_0>0$ was given in $(V_1)$. Using these numbers, we set the functions
$$
\tilde{f}(t)=
\left\{
\begin{array}{l}
f(t), \quad t \leq a, \\
\mbox{}\\
\frac{V_0}{k}t, \quad t \geq a
\end{array}
\right.
$$
and
$$
g(x,t)=\chi_\Lambda(x)f(t)+(1-\chi_\Lambda)\tilde{f}(t), \quad \forall (x,t) \in \mathbb{R}^{2},
$$
where $\Lambda$ was given in $(V_2)$ and $\chi_\Lambda$ denotes the characteristic function associated with $\Lambda$, that is, 
$$
\chi(x)=
\left\{
\begin{array}{l}
1, \quad x \in \Lambda, \\
0, \quad x \in \Lambda^c.
\end{array}
\right.
$$
Using the above functions, we will study the existence of positive solution for the following problem
$$
\left\{
\begin{array}{l}
(-\Delta)^{1/2}{u}+V(\epsilon x)u=g_\epsilon( x,u), \quad x \in \mathbb{R}, \\
\mbox{}\\
u \in H^{1/2}(\mathbb{R}),
\end{array}
\right.
\leqno{(AP)}
$$
where 
$$
g_\epsilon(x,t)=g(\epsilon x,t), \quad \forall (x,t) \in \mathbb{R}^{2}.
$$
We recall by using \cite{caffarelli}, to get a solution for the above problem, it is enough to study the existence of solution for the problem
$$
 \left\{ \begin{array}{rcl}
 -\dive ( \nabla w)=0 & \mbox{in}&  \R^{2}_{+}\noindent\\
  \frac{\partial w}{\partial \nu}=V(\epsilon x)w - g_\epsilon(x, w)  & \mbox{on}& \R
 \times\{0\},\\
  \end{array}
\right. 
\eqno{(AP)'}
$$
because if $w$ is a solution of $(AP)'$, the function $u(x)=w(x,0)$ is a solution for $(AP)$.

Here, we would like point out that if $v_\epsilon \in X^{1}(\mathbb{R}^{2}_+)$ is a solution of $(AP)'$ with
$$
v_\epsilon(x,0)<a, \quad \forall x \in \Lambda^{c}_\epsilon,
$$
where $\Lambda_\epsilon=\Lambda / \epsilon$, then $u_\epsilon(x)=v_\epsilon(x,0)$ is a solution  of \eqref{Pa}.

Associated with $(AP)'$, we have the energy functional $E_\epsilon:X^{1}(\mathbb{R}^{2}_+) \to \mathbb{R}$ given by
$$
E_\epsilon (v)=\frac{1}{2}\int_{\mathbb{R}_{+}^{2}}|\nabla v|^{2}\,\ud x \ud y+
\frac{1}{2}\int_{\mathbb{R}}V(\epsilon x)|v(x,0)|^{2}\,\ud x-\int_{\mathbb{R}}G_\epsilon(x,v(x,0))\,\ud x
$$
where
$$
G_\epsilon(x,t)=\int_{0}^{t}g_\epsilon(x,\tau)\,\ud \tau, \quad \forall (x,t) \in \mathbb{R}^{2}.
$$
Using the definition of $g$, it follows that\\

\noindent $(g_1) \quad \theta G_\epsilon(x,t) \leq g_\epsilon(x,t)t, \quad \forall (x,t) \in \Lambda_\epsilon \times \mathbb{R},$ \\
\medskip

\noindent $(g_2) \quad 2 G_\epsilon(x,t) \leq g_\epsilon(x,t)t \leq \frac{V_0}{k}|t|^{2}, \quad \forall (x,t) \in (\Lambda_\epsilon)^{c} \times \mathbb{R}.$ \\
\medskip

From assumption $(g_2)$,
$$
L(x,t)=V(x)-G_{\epsilon}(x,t) \geq \left( 1- \frac{1}{2k} \right)V(x)|t|^{2}\geq 0, \quad \forall (x,t) \in (\Lambda_\epsilon)^{c} \times \mathbb{R},
\leqno{(g_3)}
$$
$$
M(x,t)=V(x)-g_{\epsilon}(x,t)t \geq \left( 1- \frac{1}{k} \right)V(x)|t|^{2} \geq 0, \quad \forall (x,t) \in (\Lambda_\epsilon)^{c} \times \mathbb{R}.
\leqno{(g_4)}
$$

\begin{lemma} \label{L1} The functional $E_\epsilon$ verifies the mountain pass geometry, that is, \\
\noindent $i)$ \quad There are $r,\rho>0$ such that
$$
E_\epsilon(v) \geq \rho, \quad \mbox{for} \quad \|v\|=r
$$
\noindent $ii)$ \quad There is $e \in X^{1}(\mathbb{R}^{2}_+)$ with $\|e\|>r$ and $E_\epsilon(e)<0$. 
\end{lemma}
\noindent {\bf Proof.} \,  From $(g_1)-(g_4)$, there exist $c_1,c_2>0$ verifying
$$
E_\epsilon(v) \geq c_1\|v\|^{2} -c_2\|v\|^{q}, \quad \forall v \in X^{1}(\mathbb{R}^{2}_+).
$$
From the above inequality, there are $r, \rho>0$ such that
$$
E_\epsilon(v) \geq \rho, \quad \mbox{for} \quad \|v\|_{1,s}=r,
$$
showing $i)$. To prove $ii)$, fix $\varphi \in X^{1}(\R^{2}_{+})$ with $\supp{\,\varphi} \subset \Lambda_\epsilon \times \mathbb{R}$. Then, for $t>0$
$$
E_\epsilon(t \varphi)=\frac{t^{2}}{2}\|\varphi\|^{2}-\int_{\mathbb{R}}F(t\varphi(x,0))\,\ud x.
$$
From $(f_3)$, we know that there are $c_3,c_4 \geq 0$ verifying
$$
F(t) \geq c_1|t|^{\theta}-c_2, \quad \forall t \geq 0.
$$
Using the above inequality, we derive
$$
\lim_{t \to +\infty}E_\epsilon(t \varphi)=-\infty.
$$
Thereby, $ii)$ follows with $e=t \varphi$ and $t$ large enough.
\fim

In what follows, we denote by $c_\epsilon$ the mountain pass level associated with $E_\epsilon$. Related to the case $\epsilon=0$, it is possible to prove that there is $w_0 \in X^{1}(\mathbb{R}^{2}_+)$ such that 
\begin{equation} \label{exist. w}
J_0(w_0)=c_0 \quad \mbox{and} \quad J'_0(w_0)=0.
\end{equation}
The existence of $w_0$ can be obtained repeating the same approach explored in \cite{AdoOM}.
\vspace{0.5 cm}

\begin{lemma} \label{ESTIMATIVA SUPERIOR} The minimax level $c_0$ verifies  
$$
0< c_0 < \min\{1,V_0\}\left(\frac{1}{2}-\frac{1}{\theta}\right).
$$

\end{lemma}
\noindent {\bf Proof.}
Consider $w_* \in X^{1}(\mathbb{R}^{2}_+)$ verifying 
$$
\tilde{J}_{0}(w_*)=\beta_p \quad \mbox{and} \quad \tilde{J}_{0}'(w_*)=0. 
$$
By characterization of $c_0$, 
$$
c_0 \leq \max_{t \geq 0}J_0(tw_*).
$$
Consequently, by $(f_5)$, 
$$
c_0 \leq \max_{t \geq 0}\left\{ \frac{t^{2}}{2}\int_{\mathbb{R}^{2}_+}|\nabla w_*|^{2}\,\ud x \ud y+\frac{1}{2}\int_{\mathbb{R}}V_0|w_*(x,0)|^{2}\,\ud x-\frac{C_pt^{p}}{p}\int_{\mathbb{R}}|w_*(x,0)|^{p}\,\ud x\right\},
$$
which implies that
$$
c_0 \leq C_p^{{2}/{(2-p)}}\beta_p.
$$
Hence, from $(f_5)$,
$$
0<c_0< \min\{1,V_0\}\left(\frac{1}{2}-\frac{1}{\theta}\right).
$$
\fim

Hereafter, we will assume that $k$ is large enough such that
$$
0<c_0< \min\{1,V_0\}\left(\left(\frac{1}{2}-\frac{1}{\theta}\right)-\frac{1}{k}\right)<\min\{1,V_0\}\left(\frac{1}{2}-\frac{1}{\theta}\right).
$$

The next lemma establishes an important relation between $c_\epsilon$ and $c_0$.

\begin{lemma} \label{cepsilonc0} The numbers $c_0$ and $c_\epsilon$ verify the equality below
\begin{equation} \label{limite de cepsilon}
\lim_{\epsilon \to 0}c_\epsilon=c_0.
\end{equation}
Hence, there is $\epsilon_0>0$ such that 
\begin{equation} \label{limite de cepsilon2}
0<\sup_{\epsilon \in (0,\epsilon_0)}c_\epsilon < \min\{1,V_0\}\left(\left(\frac{1}{2}-\frac{1}{\theta}\right)-\frac{1}{k}\right).
\end{equation} 

\end{lemma}
\noindent {\bf Proof.} \, From $(V_1)$,
$$
c_\epsilon  \geq c_0, \quad \forall \epsilon \geq 0.
$$
Then,
\begin{equation} \label{limita1}
\liminf_{\epsilon \to 0}c_\epsilon  \geq c_0.
\end{equation}
Next, fix  $t_\epsilon >0$ such that 
$$
t_\epsilon w \in  \mathcal{M}_\epsilon=\{v \in X^{1}(\mathbb{R}^{2}_+) \setminus \{0\}\,:\, E'_\epsilon(v)v=0\}.
$$
By definition of $c_\epsilon$, we know that
$$
c_\epsilon \leq \max_{t \geq 0}E_\epsilon(tw)=E_\epsilon(t_\epsilon w).
$$
Now standard arguments as those used in \cite{OS}, it is possible to prove that
$$
\lim_{\epsilon \to 0}t_\epsilon=1
$$
and
$$
\lim_{\epsilon \to 0}E_\epsilon(t_\epsilon w)=J_0(w).
$$
Thus, 
\begin{equation} \label{limita2}
\limsup_{\epsilon \to 0}c_\epsilon \leq J_0(w)=c_0.
\end{equation}
From (\ref{limita1}) and (\ref{limita2}),
$$
\limsup_{\epsilon \to 0}c_\epsilon =c_0,
$$
showing (\ref{limite de cepsilon}). The inequality (\ref{limite de cepsilon2}) is an immediate consequence of (\ref{limite de cepsilon}) and Lemma \ref{ESTIMATIVA SUPERIOR}.
\fim

\begin{lemma} \label{Est. norma} Let $\epsilon \in (0, \epsilon_0)$ and $(v_n) \subset X^{1}(\mathbb{R}^{2}_+)$ be a $(PS)_{c_\epsilon}$sequence for $E_\epsilon$. Then,
$$
\limsup_{n \to +\infty}\|v_n\|^{2}<1.
$$
\end{lemma}

\noindent {\bf Proof.} \, Gathering $E_\epsilon(u_n)-\frac{1}{\theta}E_\epsilon'(u_n)u_n=c_\epsilon+o_n(1)$ with definition of $g$, we find 
$$ 
\left( \frac{1}{2}-\frac{1}{\theta}\right)\int_{\mathbb{R}^{2}_+}|\nabla v_n|^{2}\, \ud x\ud y+\left(\left( \frac{1}{2}-\frac{1}{\theta}\right)-\frac{1}{k} \right)V_0\int_{\mathbb{R}}|v_n(x,0)|^{2}\,\ud x \leq c_\epsilon+o_n(1),
$$
from where it follows that
$$
\min\{1,V_0\} \left(\left( \frac{1}{2}-\frac{1}{\theta}\right)-\frac{1}{k} \right) \limsup_{n \to +\infty}\|v_n\|^{2} \leq c_\epsilon < \min\{1,V_0\}\left(\left(\frac{1}{2}-\frac{1}{\theta}\right)-\frac{1}{k}\right) ,
$$
and so,
$$
\limsup_{n \to +\infty}\|v_n\|^{2} < 1.
$$
\fim

\begin{lemma} \label{L2} For $\epsilon \in (0, \epsilon_0)$, the functional $E_\epsilon$ verifies the $(PS)_{c_\epsilon}$ condition.  

\end{lemma}
\noindent {\bf Proof.} \, Let $(v_n) \subset X^{1}(\mathbb{R}^{2}_+)$ be a $(PS)_{c_\epsilon}$ sequence for $E_\epsilon$, that is,
$$
E_\epsilon(v_n) \to c_\epsilon \quad \mbox{and} \quad E'_\epsilon(v_n) \to 0, \ \mbox{as}\ n \to \infty.
$$ 
From Lemma \ref{Est. norma}, $(v_n)$ is bounded in $X^{1}(\mathbb{R}^{2}_+)$ and
$$
\limsup_{n \to +\infty}\|v_n\|^{2}<1.
$$
Since $X^{1}(\mathbb{R}^{2}_+)$ is reflexive, there is a subsequence of $(v_n)$, still denoted by  itself, and $v \in  X^{1}(\mathbb{R}^{2}_+)$ such that
$$
v_n \rightharpoonup v \quad \mbox{weakly in} \quad X^{1}(\mathbb{R}^{2}_+),  \ \mbox{as}\ n \to \infty,
$$
$$
v_n \to v \quad \mbox{in} \quad  L_{loc}^{q}(\mathbb{R}), \quad \forall q \in [2,+\infty) , \ \mbox{as}\ n \to \infty,
$$
and
$$
v_n(x,0) \to v(x,0) \quad \mbox{a.e. in} \quad \mathbb{R}, \ \mbox{as}\ n \to \infty.
$$
Moreover, by Lemma \ref{Convergencia em limitados},
$$
\int_{\mathbb{R}}f(v_n(x,0))\phi(x,0)\,\ud x \to \int_{\mathbb{R}}f(v(x,0))\phi(x,0)\,\ud x,
$$
as $ n \to \infty,$ for all $\phi \in C^{\infty}_{0}(\overline{\R^{2}_{+}})$.

Using the above limits, it is possible to prove that $v$ is a critical point for $E_\epsilon$, that is, 
$$
E'_\epsilon(v)\varphi=0, \quad \forall \varphi \in X^{1}(\mathbb{R}^{2}_+).
$$
Considering $\varphi=v$, we have that $E'_\epsilon(v)v=0$, and so,
\begin{eqnarray*}
&&\int_{\mathbb{R}_{+}^{2}}|\nabla v|^{2}\,\ud x \ud y+\int_{\Lambda_\epsilon}V(\epsilon x)|v(x,0)|^{2}\,\ud x+
\int_{(\Lambda_\epsilon)^{c}}M(x,v(x,0))\,\ud x  \\ &&=\int_{\Lambda_\epsilon}f(v(x,0))v(x,0)\,\ud x.
\end{eqnarray*}
On the other hand, using the limit $E'_\epsilon(v_n)v_n=o_n(1)$, we derive that
\begin{eqnarray*}
&&\int_{\mathbb{R}_{+}^{2}}|\nabla v_n|^{2}\,\ud x\ud y+\int_{\Lambda_\epsilon}V(\epsilon x)|v_n(x,0)|^{2}\,\ud x+
\int_{(\Lambda_\epsilon)^{c}}M(x,v_n(x,0))\,\ud x \\ && =\int_{\Lambda_\epsilon}f(v_n(x,0))v_n(x,0)\,\ud x+o_n(1).
\end{eqnarray*}
Since $\Lambda_\epsilon$ is bounded, by the compactness of sobolev embedding and Lemma \ref{Convergencia em limitados} yield
$$
\lim_{n \to +\infty}\int_{\Lambda_\epsilon}f(v_n(x,0))v_n(x,0)\,\ud x=\int_{\Lambda_\epsilon}f(v(x,0))v(x,0)\,\ud x
$$
and
\begin{equation} \label{Eq1}
\lim_{n \to +\infty}\int_{\Lambda_\epsilon}V(\epsilon x)|v_n(x,0)|^{2}\,\ud x=\int_{\Lambda_\epsilon}V(\epsilon x)|v(x,0)|^{2}\,\ud x.
\end{equation}
Therefore,
\begin{eqnarray*}
&&\limsup_{n \to +\infty}\left( \int_{\mathbb{R}_{+}^{2}}|\nabla v_n|^{2}\,\ud x\ud y+\int_{(\Lambda_\epsilon)^{c}}M(x,v_n(x,0))
\,\ud x \right)\\ && =\int_{\mathbb{R}_{+}^{2}}|\nabla v|^{2}\,\ud x\ud y+\int_{(\Lambda_\epsilon)^{c}}M(x,v(x,0))\,\ud x.
\end{eqnarray*}
Now, recalling that $M(x,t) \geq 0$, the Fatous' lemma leads
\begin{eqnarray*}
&&\
\liminf_{n \to +\infty}\left( \int_{\mathbb{R}_{+}^{2}}|\nabla v_n|^{2}\,\ud x\ud y+\int_{(\Lambda_\epsilon)^{c}}M(x,v_n(x,0))
\,\ud x \right)\\&& \geq \int_{\mathbb{R}_{+}^{2}}|\nabla v|^{2}\,\ud x\ud y+\int_{(\Lambda_\epsilon)^{c}}M(x,v(x,0))\,\ud x
\end{eqnarray*}
Hence,
\begin{equation} \label{Eq2}
\lim_{n \to +\infty}\int_{\mathbb{R}_{+}^{2}}|\nabla v_n|^{2}\,\ud x\ud y=\int_{\mathbb{R}_{+}^{2}}|\nabla v|^{2}\,\ud x\ud y
\end{equation}
and
$$
\lim_{n \to +\infty}\int_{(\Lambda_\epsilon)^{c}}M(x,v(x,0))\,\ud x=\int_{(\Lambda_\epsilon)^{c}}M(x,v(x,0))\,\ud x.
$$
The last limit combined with definition of function $M$ gives
$$
\lim_{n \to +\infty}\int_{(\Lambda_\epsilon)^{c}}V(\epsilon x)|v_n(x,0)|^{2}\,\ud x=\int_{(\Lambda_\epsilon)^{c}}
V(\epsilon x)|v(x,0)|^{2}\,\ud x.
$$
Gathering this limit with (\ref{Eq1}), we deduce that
\begin{equation} \label{Eq3}
\lim_{n \to +\infty}\int_{\mathbb{R}}V(\epsilon x)|v_n(x,0)|^{2}\,\ud x=\int_{\mathbb{R}}V(\epsilon x)|v(x,0)|^{2}\,\ud x.
\end{equation} 
From (\ref{Eq2})-(\ref{Eq3}),
$$
\lim_{n \to +\infty}\|v_n\|_{\epsilon}^{2}=\|v\|_{\epsilon}^{2}.
$$
As $X^{1}(\mathbb{R}^{2}_+)$ is a Hilbert space and $v_n \rightharpoonup v$ weakly in $X^{1}(\mathbb{R}^{2}_+)$, as $n \to \infty,$ the above limit yields
$$
v_n \to v \quad \mbox{in} \quad  X^{1}(\mathbb{R}^{2}_+), \ \mbox{as}\ n \to \infty,
$$
showing that $E_\epsilon$ verifies the $(PS)_{c_\epsilon}$.
\fim

\begin{theorem} \label{T0} For $\epsilon \in (0,\epsilon_0)$, the functional $E_\epsilon$ has a nonnegative critical point $v_\epsilon \in X^{1}(\R^2_{+})$ such 
\begin{equation} \label{vepsilon}
E_\epsilon(v_\epsilon)=c_\epsilon \quad \mbox{and} \quad E'_\epsilon(v_\epsilon)=0.
\end{equation}
\end{theorem}
\noindent {\bf Proof.} \, From Lemma \ref{cepsilonc0}, there is $\epsilon_0>0$, such that $E_\epsilon$ verifies the $(PS)_{c_\epsilon}$ condition for $\epsilon \in (0, \epsilon_0).$ Then, the existence of  $v_\epsilon $ is an immediate consequence of  the Mountain Pass Theorem due to Ambrosetti and Rabinowitz (see e.g.  \cite{willem}). The function $v_\epsilon$ is nonnegative, because
$$
E'_\epsilon(v_\epsilon)(v_{\epsilon}^{-})=0  \Longrightarrow v_{\epsilon}^{-}=0,
$$
where $v_{\epsilon}^{-}=\min\{v_{\epsilon},0\}$.\fim

\begin{lemma} \label{sequencia} Decreasing $\epsilon_0$, if necessary, there are $r, \beta>0$ and $(y_\epsilon) \subset \mathbb{R}$ such that
\begin{equation}\label{Ciri}
\int_{y_\epsilon-r}^{y_\epsilon+r}|v_\epsilon(x,0)|^{2}\,\ud x\geq \beta, \quad \forall \epsilon \in (0, \epsilon_0). 
\end{equation}
\end{lemma}
\noindent {\bf Proof.} \,  First of all, we recall that since $(v_\epsilon)$ satisfies $(\ref{vepsilon}),$
 there is $\alpha >0$, which is independent of $\epsilon$, such that
\begin{equation} \label{alfa}
\|v_\epsilon\|_{\epsilon}^{2} \geq \alpha, \quad \forall \epsilon >0.
\end{equation}

To show \eqref{Ciri}, it is enough to see that for any sequence $(\epsilon_n) \subset (0,+\infty)$ with $\epsilon_n \to 0$, the limit below
$$
\lim_{n \to +\infty}\sup_{y \in \mathbb{R}}\int_{y-r}^{y+r}|v_{\epsilon_n}(x,0)|^{2}\,\ud x=0,
$$
does not hold for any $r>0$. Otherwise, if it holds for some $r>0$, by Lemma \ref{lions},
$$
\int_{\R}f(v_{\epsilon_n}(.,0))v_n(x,0) \ud x \to 0, \ \mbox{as}\  n \to \infty,
$$
implying that 
$$
\|v_{\epsilon_n}\|_{\epsilon}^{2} \to 0 \quad \mbox{as} \quad n \to +\infty,
$$
which contradicts (\ref{alfa}). 

\fim

\begin{lemma} \label{convergencia em Lambda} For any $\epsilon_n \to 0$, consider the sequence $(y_{\epsilon_n}) 
\subset \mathbb{R}$ given in Lemma \ref{sequencia} and $\psi_n(x,y)=v_{\epsilon_n}(x+y_{\epsilon_n},y).$ Then, up to subsequence,  there is $\psi \in X^{1}(\mathbb{R}^{2}_+$ such that
\begin{equation} \label{psi}
\psi_n \to \psi \quad \mbox{in} \quad  X^{1}(\mathbb{R}^{2}_+), \ \mbox{as}\  n \to \infty.
\end{equation}
Moreover, there is $x_0 \in \Lambda$ such that
\begin{equation} \label{yn}
\lim_{n \to 0}\epsilon_n y_{\epsilon_n}=x_0 \quad \mbox{and} \quad V(x_0)=V_0.
\end{equation} 
\end{lemma}
\noindent {\bf Proof.} \, We begin the proof showing that  $(\epsilon_n y_{\epsilon_n})$ is a bounded sequence. Hereafter, we denote by $(y_n)$ and $(v_n)$ the sequences $(y_{\epsilon_n})$ and $(v_{\epsilon_n})$ respectively. 

Since $E'_{\epsilon_n}(v_n)\phi=0, \forall \phi \in X^{1}(\mathbb{R}^{2}_+)$, we have that
$$
\int_{\mathbb{R}_{+}^{2}}\nabla v_n \nabla \phi \,\ud x\ud y+\int_{\mathbb{R}}V(\epsilon_n x)v_n(x,0) \phi(x,0)\,\ud x-
\int_{\mathbb{R}}g_\epsilon(x,v_n(x,0))\phi(x,0)\,\ud x=0.
$$
Then, 
$$
\int_{\mathbb{R}_{+}^{2}}|\nabla v_n|^{2}\,\ud x\ud y+\int_{\mathbb{R}}V(\epsilon_n x)|v_n(x,0)|^{2}\,\ud x-
\int_{\mathbb{R}}g_\epsilon(x,v_n(x,0))v_n(x,0)\,\ud x=0.
$$
From definition of $g$, we see that
$$
g_\epsilon(x,t)\leq f(t), \quad \forall t \geq 0,
$$
and reminding that $v_n \geq 0,$ we infer that 
$$
\int_{\mathbb{R}_{+}^{2}}|\nabla v_n|^{2}\,\ud x\ud y+\int_{\mathbb{R}}V_0|v_n(x,0)|^{2}\,\ud x-
\int_{\mathbb{R}}f(v_n(x,0))v_n(x,0)\,\ud x\leq0.
$$
Therefore, there is $s_n \in (0,1)$ such that
$$
s_nv_n \in \mathcal{M}_0=\{v \in X^{1}(\mathbb{R}^{2}_+)\setminus \{0\}\,:\, J'_0(v)v=0\}.
$$
Using the characterization of $c_0$, we know that
$$
c_0 \leq J_0(s_nv_n), \quad \forall n \in \mathbb{N}.
$$
As 
$$
J_0(w) \leq E_\epsilon(w), \quad \forall w \in X^{1}(\mathbb{R}^{2}_+) \quad \mbox{and} \quad \epsilon >0,
$$
it follows that
$$
c_0 \leq J_0(s_nv_n) \leq E_{\epsilon_n}(s_nv_n)\leq \max_{s \geq 0}E_{\epsilon_n}(sv_n)=E_{\epsilon_n}(v_n)=c_{\epsilon_n}. 
$$
Recalling that
$$
c_{\epsilon_n} \to c_0,\ \mbox{as}\  n \to \infty,
$$
the last inequality gives
$$
(s_nv_n) \subset \mathcal{M}_0, \ \forall n \in \N,  \quad \mbox{and} \quad J_0(s_nv_n) \to c_0, \ \mbox{as}\  n \to \infty.
$$
By change variable, we also have
$$
(s_n\psi_n) \subset \mathcal{M}_0, \ \forall n \in \N, \quad \mbox{and} \quad J_0(s_n\psi_n) \to c_0, \ \mbox{as}\  n \to \infty.
$$
Using Ekeland Variational Principle, we can assume that $(s_nv_n)$ is a $(PS)_{c_0}$ sequence, that is,
$$
(s_n\psi_n) \subset \mathcal{M}_0, \forall n \in \N, \quad J_0(s_n\psi_n) \to c_0 \quad \mbox{and} \quad J'_0(s_n\psi_n) \to 0, \ \mbox{as}\  n \to \infty.
$$
A direct computation shows that $(s_n)$ is a bounded sequence with
$$
\liminf_{n \to +\infty}s_n>0.
$$
Thus, in what follows, we can assume that for some subsequence, there is $s_0>0$ such that
$$
s_n \to s_0, \ \mbox{as}\  n \to \infty.
$$
From definition of $y_n$ and $\psi_n$, we know that $\psi \in X^{1}(\mathbb{R}^{2}_+) \setminus \{0\}$. Moreover, as $J'_0(s_n\psi_n) \to 0$, 
we also have $ J'_0(s_0\psi)=0$. Thereby, by definition of $c_0$, we obtain
$$
c_0 \leq J_0(s_0\psi).
$$ 
On the other hand, by Fatou's Lemma we obtain
$$
\liminf_{n \to +\infty}J_0(s_n\psi_n) \geq J_0(s_0\psi)
$$
which implies
$$
J_0(s_0\psi)=c_0 \quad \mbox{and} \quad J'_0(s_0\psi)=0.
$$
The above equalities combined with Fatous' Lemma, up to a subsequence,  gives 
$$
s_n \psi_n \to s_0 \psi \quad \mbox{in} \quad X^{1}(\mathbb{R}^{2}_+), \ \mbox{as}\  n \to \infty.
$$
Recalling that $s_n \to s_0>0, \ \mbox{as}\  n \to \infty$, we can conclude that
$$
\psi_n \to \psi \quad \mbox{in} \quad X^{1}(\mathbb{R}^{2}_+), \ \mbox{as}\  n \to \infty,
$$ 
showing (\ref{psi}). 

Using the last limit, we are able to prove (\ref{yn}). To do end, we begin making the following claim

\begin{claim} \label{convergencia yn} $\displaystyle \lim_{n \to +\infty}\dist(\epsilon_n y_n, \overline{\Lambda})=0$

\end{claim}

Indeed, if the claim does not hold, there is $\delta>0$ and a subsequence of $(\epsilon_n y_n)$, still denoted by itself, such that,
$$
\dist(\epsilon_n y_n, \overline{\Lambda}) \geq \delta, \quad \forall n \in \mathbb{N}.
$$
Consequently, there is $r>0$ such that
$$
(\epsilon_n y_n-r,\epsilon_n y_n+r) \subset \Lambda^{c}, \quad \forall n \in \mathbb{N}.
$$ 
From definition of $\psi_n$, we have that 
\begin{eqnarray*}
&&\int_{\mathbb{R}_{+}^{2}}|\nabla \psi_n|^{2}\,\ud x\ud y+\int_{\mathbb{R}}V(\epsilon_n x +
\epsilon_n y_n)|\psi_n(x,0)|^{2}\,\ud x \\ && =\int_{\mathbb{R}}g(\epsilon_nx +\epsilon_n y_n,\psi_n(x,0))\psi_n(x,0)\,\ud x.
\end{eqnarray*}
Note that
\begin{eqnarray*}
&&
\int_{\mathbb{R}}g(\epsilon_nx +\epsilon_n y_n,\psi_n(x,0))\psi_n(x,0)\,\ud x  \leq 
\int_{-\frac{r}{\epsilon_n}}^{\frac{r}{\epsilon_n}} g(\epsilon_nx +\epsilon_n y_n,\psi_n(x,0))\psi_n(x,0)\,\ud x \\ && +
(\int_{-\infty}^{-\frac{r}{\epsilon_n}}+\int_{\frac{r}{\epsilon_n}}^{+\infty} )g(\epsilon_nx +\epsilon_n y_n,\psi_n(x,0))\psi_n(x,0)\,\ud x,
\end{eqnarray*}
and  so,
\begin{eqnarray*}
&&
\int_{\mathbb{R}}g(\epsilon_nx +\epsilon_n y_n,\psi_n(x,0))\psi_n(x,0)\,\ud x \\ &&  \leq \frac{V_0}{k}
\int_{-\frac{r}{\epsilon_n}}^{\frac{r}{\epsilon_n}}|\psi_n(x,0)|^{2}\,\ud x+(\int_{-\infty}^{-\frac{r}{\epsilon_n}}+\int_{\frac{r}{\epsilon_n}}^{+\infty} )
f(\psi_n(x,0))\psi_n(x,0)\,\ud x.
\end{eqnarray*}
Therefore,
\begin{eqnarray*}
&&
\int_{\mathbb{R}_{+}^{2}}|\nabla \psi_n|^{2}\,\ud x\ud y+\int_{\mathbb{R}}V(\epsilon_n x +\epsilon_n y_n)|\psi_n(x,0)|^{2}
\,\ud x\\ &&  \leq 
\frac{V_0}{k}\int_{-\frac{r}{\epsilon_n}}^{\frac{r}{\epsilon_n}}|\psi_n(x,0)|^{2}\,\ud x+(\int_{-\infty}^{-\frac{r}{\epsilon_n}}+\int_{\frac{r}{\epsilon_n}}^{+\infty} )f(\psi_n(x,0))\psi_n(x,0)\,\ud x.
\end{eqnarray*}
implying that
\begin{equation} \label{Eq4}
\int_{\mathbb{R}_{+}^{2}}|\nabla \psi_n|^{2}\,\ud x\ud y+A \int_{\mathbb{R}^{N}}|\psi_n(x,0)|^{2}\,\ud x \leq
(\int_{-\infty}^{-\frac{r}{\epsilon_n}}+\int_{\frac{r}{\epsilon_n}}^{+\infty} )f(\psi_n(x,0))\psi_n(x,0)\,\ud x,
\end{equation}
where $A=V_0\left(1-\frac{1}{k}\right)$. By (\ref{psi}), 
$$
(\int_{-\infty}^{-\frac{r}{\epsilon_n}}+\int_{\frac{r}{\epsilon_n}}^{+\infty} )f(\psi_n(x,0))\psi_n(x,0)\,\ud x \to 0, \ \mbox{as} \ n \to \infty,
$$
and $ \ \mbox{as} \ n \to \infty,$
\[
\int_{\mathbb{R}_{+}^{2}}|\nabla \psi_n|^{2}\,\ud x\ud y+A \int_{\mathbb{R}}|\psi_n(x,0)|^{2}\,\ud x   \to 
\int_{\mathbb{R}_{+}^{2}}|\nabla \psi|^{2}\,\ud x\ud y+A\int_{\mathbb{R}}|\psi(x,0)|^{2}\,\ud x>0,
\]
which contradicts (\ref{Eq4}). This proves Claim \ref{convergencia yn}. 

From Claim \ref{convergencia yn}, there are a subsequence of $(\epsilon_n y_n)$ and $x_0 \in \overline{\Lambda}$ such that
$$
\lim_{n \to +\infty}\epsilon_n y_n=x_0. 
$$
\begin{claim} \label{convergencia yn 2}  $x_0 \in \Lambda$. 
\end{claim}

Indeed, from definition of $\psi_n$,
$$
\int_{\mathbb{R}_{+}^{2}}|\nabla \psi_n|^{2}\,\ud x \ud y+
\int_{\mathbb{R}}V(\epsilon_n x +\epsilon_n y_n)|\psi_n(x,0)|^{2}\,\ud x\leq \int_{\mathbb{R}}f(\psi_n(x,0))\psi_n(x,0)\,\ud x.
$$
Then, by (\ref{psi}),
$$
\int_{\mathbb{R}_{+}^{2}}|\nabla \psi|^{2}\,\ud x\ud y+\int_{\mathbb{R}}V(x_0)|\psi(x,0)|^{2}\,\ud x\leq 
\int_{\mathbb{R}}f(\psi(x,0))\psi(x,0)\,\ud x.
$$
Hence, there is $s_1 \in (0,1)$ such that
$$
s_1 \psi \in \mathcal{M}_{V(x_0)}=\left\{v \in X^{1}(\mathbb{R}^{2}_+)\setminus \{0\}\,:\,\tilde{J}'_{V(x_0)}v=v \right\}
$$
where $\tilde{J}_{V(x_0)}:X^{1}(\mathbb{R}^{2}_+) \to \mathbb{R}$ is given by 
$$
\tilde{J}_{V(x_0)} (v)=\frac{1}{2}\int_{\mathbb{R}_{+}^{2}}|\nabla v|^{2}\,\ud x\ud y+
\frac{1}{2}\int_{\mathbb{R}}V(x_0)|v(x,0)|^{2}\,\ud x-\int_{\mathbb{R}}F(v(x,0))\,\ud x.
$$
If $\tilde{c}_{V(x_0)}$ denotes the mountain pass level associated with $\tilde{J}_{V(x_0)}$, we must have 
$$
\tilde{c}_{V(x_0)} \leq \tilde{J}_{V(x_0)}(s_1\psi)\leq \liminf_{n \to +\infty}E_{\epsilon_n}(v_n)=\liminf_{n \to +\infty}c_{\epsilon_n}=c_0=\tilde{c}_{V(0)}.
$$
Hence,
$$
\tilde{c}_{V(x_0)} \leq \tilde{c}_{V(0)}, 
$$
from where it follows that
$$
V(x_0) \leq V(0).
$$
As $V_0=\inf_{x \in \mathbb{R}}V(x)$, the above inequality implies that
$$
V(x_0)=V(0)=V_0.
$$
Moreover, by $(V_2)$, $x_0 \notin \partial \Lambda$. Then, $x_0 \in \Lambda$, finishing the proof.
\fim

\begin{cor} \label{Cor Limita} Let $(\psi_n)$ the sequence given in Lemma \ref{convergencia em Lambda}. Then, $\psi_n(\cdot,0) \in L^{\infty}(\mathbb{R})$ and there is $K>0$ such that
\begin{equation} \label{Linfinito}
|\psi_n(\cdot,0)|_{\infty} \leq K, \quad \forall n \in \mathbb{N}
\end{equation}
and
\begin{equation} \label{Linfinito2}
\psi_n(\cdot,0) \to \psi(\cdot,0) \quad \mbox{in} \quad L^{p}(\mathbb{R}), \quad \forall p \in (2,+\infty), \ \mbox{as} \ n \to \infty.
\end{equation}
As an immediate consequence, the sequence $h_n(x)= g(\epsilon_n x+\epsilon_n y_n,\psi_n(x,0))$ must verify
\begin{equation} \label{Linfinito3}
h_n \to f(\psi(\cdot,0)) \quad \mbox{in} \quad L^{p}(\mathbb{R}), \quad \forall p \in (2,+\infty), \ \mbox{as} \ n \to \infty.
\end{equation}
\end{cor}

\noindent {\bf Proof. }\, In what follows,  for each $L>0$, we set  
$$
\psi_{n,L}(x,y)=
\left\{
\begin{array}{lcr} 
\psi_n(x,y), &\mbox{if}& \psi_n(x,y)\leq L\\
L, &\mbox{if}& \psi_n(x,y) \geq L
\end{array}
\right.
$$
and
$$
z_{n,L}=\psi_{n,L}^{2(\beta -1)}\psi_n,
$$
with $\beta >1$ to be determined later. Since 
\begin{eqnarray*} &&
\int_{\mathbb{R}_{+}^{2}}\nabla \psi_n \nabla \phi \,\ud x\ud y+\int_{\mathbb{R}}V(\epsilon_n x +
 \epsilon_n y_n)\psi_n(x,0) \phi(x,0)\,\ud x \\ && -\int_{\mathbb{R}}g(\epsilon_n x + \epsilon_n y_n,\psi_n(x,0))\phi(x,0)\,\ud x=0, \quad \forall \phi \in X^{1}(\mathbb{R}^{2}_+), \ \forall n \in \N,
\end{eqnarray*}
adapting the same approach explored in C.~Alves and G.~Figueiredo \cite[Lemma 4.1]{AF2} and using the fact that $(\psi_n)$ is bounded in $X^{1}(\mathbb{R}^{2}_+)$, we conclude that there is $K>0$ such that
$$
|\psi_n(.,0)|_{\infty} \leq K, \quad \forall n \in \mathbb{N}.
$$

Now, the limit (\ref{Linfinito2}) is obtained  by interpolation on the $L^{p}$ spaces, while that (\ref{Linfinito3}) follows 
combining the growth  condition on  $g$ with (\ref{Linfinito2}). 
\fim

In what follows, we denote by $(w_n) \subset H^{1/2}(\mathbb{R})$ the sequence $(\psi_n(\cdot,0))$, that is, 
$$
w_n(x)=\psi_n(x,0), \quad \forall x \in \mathbb{R}.
$$
Since 
\begin{eqnarray*} &&
\int_{\mathbb{R}_{+}^{2}}\nabla \psi_n \nabla \phi \,\ud x\ud y+\int_{\mathbb{R}}V(\epsilon_n x +
 \epsilon_n y_n)\psi_n (x,0) \phi(x,0)\,\ud x \\ && -\int_{\mathbb{R}}g(\epsilon_n x + \epsilon_n y_n,\psi_n(x,0))\phi(x,0)\,\ud x=0, \quad \forall \phi \in X^{1}(\mathbb{R}^{2}_+),
\end{eqnarray*}
we have that $w_n$ is a solution of the problem
$$
(-\Delta)^{1/2}{w_n}+V(\epsilon_n x+\epsilon_n y_n)w_n=g(\epsilon_n x+\epsilon_n y_n,w_n), \quad \mbox{in} \quad \mathbb{R},
$$
or equivalently,
\begin{equation} \label{C1}
(-\Delta)^{1/2}{u}+w_n=\chi_n, \quad \mbox{in} \quad \mathbb{R},
\end{equation}
where
\begin{equation} \label{C3}
\chi_n(x)=w_n(x)+g(\epsilon_n x+\epsilon_n y_n x,w_n(x))-V(\epsilon_n x+\epsilon_n y_n)w_n(x), \quad x \in \mathbb{R}.
\end{equation}
Denoting $\chi(x)=w(x)+f(w(x))-V(x_0)w(x)$, by Corollary \ref{Cor Limita}, we have that
\begin{equation} \label{C3}
\chi_n \to \chi \quad \mbox{in} \quad L^{p}(\mathbb{R}), \quad \forall p \in [2,+\infty), \ \mbox{as} \ n \to \infty,
\end{equation}
and there is $k_1>0$,
\begin{equation} \label{C4}
|\chi_n|_{\infty} \leq k_1, \quad \forall n \in \mathbb{N}.
\end{equation}

Motivated by some results found in \cite{BGe} ( see also \cite{FQT} ),  which holds for whole line, we deduce that
$$
w_n(x)=(\mathcal{K}\ast \chi_n)(x)=\int_{\mathbb{R}}\mathcal{K}(x-y)\chi_n(y)\; \ud y,
$$
where $\mathcal{K}$ is the Bessel kernel, which verifies:
\begin{flushleft}
\noindent $(K_1)$ \; $\mathcal{K} $ is positive and even on $\mathbb{R} \setminus \{0\}$, \\
\noindent $(K_2)$ \; There is $C>0$ such that  
$$
\mathcal{K}(x) \leq C/|x|^{2}, \quad \forall x \in \mathbb{R} \setminus \{0\} 
$$
\noindent and \\
$(K_3)$ \; $\mathcal{K} \in L^{q}(\mathbb{R}), \quad \forall q \in L^{q}(\mathbb{R}) \quad \forall q \in [1, \infty]$. \\
\end{flushleft}

Using the above informations, we are able to prove the following result

\begin{lemma} \label{Con. ZERO}The sequence $(w_n)$ verifies 
$$	
w_n(x) \to 0 \quad \mbox{as} \quad |x| \to +\infty,
$$
uniformly in $n \in \mathbb{N}.$
\end{lemma}
\noindent {\bf Proof.} Given $\delta >0$, we have
\begin{eqnarray*}
0\leq w_n(x) &\leq&  \int_{\mathbb{R}}\mathcal{K}(x-y)|\chi_n|(y)\; \ud y\\ & =& (\int_{-\infty}^{x-{1}/{\delta}}+\int_{x+{1}/{\delta}}^{+\infty})\mathcal{K}(x-y)|\chi_n|(y)\;\ud y+\int_{x-{1}/{\delta}}^{x+{1}/{\delta}}\mathcal{K}(x-y)|\chi_n|(y)\;\ud y
\end{eqnarray*}
from  $(K_2)$, we have that, for all $n \in \N,$
\begin{eqnarray} \label{C5}
 (\int_{-\infty}^{x-{1}/{\delta}}+\int_{x+{1}/{\delta}}^{+\infty})\mathcal{K}(x-y)|\chi_n|(y)\;\ud y &\leq& C\delta^{1/2}|\chi_n|_\infty (\int_{-\infty}^{x-{1}/{\delta}}+\int_{x+{1}/{\delta}}^{+\infty})\frac{\;\ud y}{|x-y|^{3/2}} \\
&\leq& C\delta^{1/2}k_1 (\int_{-\infty}^{x-1}+\int_{x+1}^{+\infty})\frac{\;}{|x-y|^{3/2}}=C_1\delta^{1/2}.\nonumber
\end{eqnarray}
On the other hand,
$$
\int_{x-{1}/{\delta}}^{x+{1}/{\delta}}\mathcal{K}(x-y)|\chi_n|(y)\;\ud y\leq \int_{x-{1}/{\delta}}^{x+{1}/{\delta}}\mathcal{K}(x-y)|\chi_n - \chi|(y)\;\ud y+\int_{x-{1}/{\delta}}^{x+{1}/{\delta}}\mathcal{K}(x-y)|\chi|(y)\;\ud y.
$$
Fix $q >1$ with $q$ sufficiently close to $ 1$ and $q'>2$ such that $1/q+1/q'=1$. From $(K_2)$ and (\ref{C1}), 
$$
\int_{x-{1}/{\delta}}^{x+{1}/{\delta}}\mathcal{K}(x-y)|\chi_n|(y)\;\ud y\leq |\mathcal{K}|_{q}|\chi_n-\chi|_{q'}+|\mathcal{K}|_{q}|\chi|_{L^{q'}(x-{1}/{\delta},x+{1}/{\delta})}
$$ 
As 
$$
|\chi_n-\chi|_{q'} \to 0 \quad \mbox{as} \quad n \to +\infty
$$
and
$$
|\chi|_{L^{q'}(x-{1}/{\delta},\, x+{1}/{\delta})} \to 0 \quad \mbox{as} \quad |x| \to +\infty,
$$
we deduce that there are  $R>0$ and $n_0 \in \mathbb{N}$ such that
\begin{equation} \label{C6}
\int_{x-{1}/{\delta}}^{x+{1}/{\delta}}\mathcal{K}(x-y)|\chi_n|(y)\;\ud y\leq \delta, \quad \forall  n \geq n_0 \quad \mbox{and} \quad |x|\geq R.
\end{equation}
from (\ref{C5}) and (\ref{C6}),
\begin{equation} \label{C7}
\int_{\mathbb{R}}\mathcal{K}(x-y)|\chi_n|(y)\;\ud y\leq C_1\delta^{d}+\delta, \quad \forall  n \geq n_0 \quad \mbox{and} \quad |x|\geq R.
\end{equation}
The same approach can be used to prove that for each $n \in \{1,...., n_{0}-1\}$, there is $R_n>0$ such that
\begin{equation} \label{C8}
\int_{\mathbb{R}}\mathcal{K}(x-y)|\chi_n|(y)\;\ud y\leq C_1\delta^{d}+\delta, \quad  |x|\geq R_n.
\end{equation}
Hence, increasing $R, $ if necessary, we must have
$$
\int_{\mathbb{R}}\mathcal{K}(x-y)|\chi_n|(y)\;\ud y\leq C_1\delta^{d}+\delta, \quad \mbox{for} \quad |x|\geq R, \quad \mbox{uniformly in } \quad n \in \mathbb{N}.
$$	
Since $\delta$ is arbitrary, the proof is finished.  \fim

\vspace{0.5 cm}

\begin{cor} \label{Cor original} There is $n_0 \in \mathbb{N}$ such that
$$
v_n(x,0) < a, \quad \forall n \geq n_0 \quad \mbox{and} \quad x \in \Lambda_{\epsilon_n}^{c}.
$$	
Hence, $u_n(x)=v_n(x,0)$ is a solution of $(P'_{\epsilon_n})$ for $n \geq n_0$.
\end{cor}
\noindent {\bf Proof.}\, By Lemma \ref{convergencia em Lambda}, we know that $\epsilon_n y_n \to x_0 
$, for some $x_0 \in \Lambda$. Thereby, there is $r>0$ such that some subsequence, still denoted by itself,
$$
(r-\epsilon_n y_n,\, r+\epsilon_n y_n) \subset \Lambda, \quad \forall n \in \mathbb{N}.
$$ 
Hence, 
$$
(y_n- {r}/{\epsilon_n},\, y_n+{r}/{\epsilon_n})\subset \Lambda_{\epsilon_n}, \quad \forall n \in \mathbb{N},
$$
or equivalently
$$
\Lambda_{\epsilon_n}^{c} \subset (-\infty, y_n-{r}/{\epsilon_n})\cup (y_n+{r}/{\epsilon_n},+\infty), \quad \forall n \in \mathbb{N}.
$$
Now, by Lemma \ref{Con. ZERO}, there is $R>0$ such that
$$
w_n(x) < a, \quad \mbox{for} \quad |x| \geq R \quad \mbox{and} \quad \forall n \in \mathbb{N},
$$
from where it follows,
$$
v_n(x,0)=\psi_n(x-y_n,0)=w_n(x-y_n)<a, \quad \mbox{for} \quad x \in  (-\infty, y_n-R)\cup (y_n+R,+\infty)
$$
 and \quad $\forall n \in \mathbb{N}.$

On the other hand, we have that
$$
\Lambda^{c}_{\epsilon_n} \subset(-\infty,\, y_n-{r}/{\epsilon_n})\cup (y_n+{r}/{\epsilon_n},+\infty), \quad \forall n \in \mathbb{N}.
$$
Thus, there is $n_0 \in \mathbb{N}$, such that
$$
(-\infty, y_n-{r}/{\epsilon_n})\cup (y_n+{r}/{\epsilon_n},+\infty) \subset (-\infty, y_n-R)\cup (y_n+R,+\infty), \quad \forall n \geq n_0,
$$
implying that
$$
v_n(x,0)<a, \quad \forall x \in \Lambda^{c}_{\epsilon_n}  \quad \mbox{and} \quad n \geq n_0,
$$
finishing the proof. 
\fim

%-----------------------------------------------------------------------------------
%-----------------------------------------------------------------------------------
\section{Proof of Theorem \ref{resultado}} \label{proof-Thm}
%-----------------------------------------------------------------------------------
%-----------------------------------------------------------------------------------

By Theorem \ref{T0}, we know that problem $(AP)$ has a nonnegative solution $v_\epsilon$ for all $\epsilon>0$.   Applying Corollary \ref{Cor original}, there is $\epsilon_0$ such that
$$
v_{\epsilon}(x,0)<a, \quad \forall x \in \Lambda^{c}_{\epsilon} \quad \mbox{and} \quad \forall \epsilon \in (0, \epsilon_0),
$$
that is, $v_{\epsilon}(\cdot,0)$ is a solution of \eqref{Pa} for $\epsilon  \in (0, \epsilon_0)$. Considering
$$
u_\epsilon(x)=v_{\epsilon}({x}/{\epsilon},0), \quad \mbox{for} \quad \forall \epsilon  \in (0, \epsilon_0),
$$
is a solution for original problem  \eqref{Pb}.

If $x_\epsilon$ denotes a global maximum point of $u_\epsilon$, it is easy to see that there is $\tau_0>0$ such that 
$$
u_\epsilon(x_\epsilon) \geq \tau_0, \quad \forall \epsilon >0.
$$
In what follows, setting  $z_\epsilon= {(x_\epsilon -\epsilon y_\epsilon)}{\epsilon^{-1}}$, we have that $z_\epsilon$ is a global maximum point of $w_\epsilon$ and  
$$
w_\epsilon(z_\epsilon) \geq \tau_0, \quad \forall \epsilon >0.
$$

Now, we claim that
\begin{equation} \label{DV}
\lim_{\epsilon \to 0}V(x_\epsilon)=V_0.
\end{equation}
Indeed, by Lemma \ref{Con. ZERO}, we know that
$$
w_{\epsilon_n}(x) \to 0 \quad \mbox{as} \quad |x| \to +\infty \quad \mbox{uniformly in} \quad n \in \mathbb{N}.
$$
Therefore, $(z_\epsilon)$ is a bounded sequence. Moreover, for some subsequence, we also know that there is $x_0 \in \Lambda$ satisfying $V(x_0)=V_0$ and 
$$
\epsilon_n y_{\epsilon_n } \to x_0, \ \mbox{as} \ n \to \infty.
$$
Hence,
$$
x_{\epsilon_n}=\epsilon_n z_{\epsilon_n}+\epsilon_n y_{\epsilon_n } \to x_0, \ \mbox{as} \ n \to \infty,
$$
implying that
$$
V(x_{\epsilon_n}) \to V_0, \ \mbox{as} \ n \to \infty,
$$
showing that (\ref{DV}) holds.

%-----------------------------------------------------------------------------------
%-----------------------------------------------------------------------------------

\bigskip
%-----------------------------------------------------------------------------------
%-----------------------------------------------------------------------------------


\begin{thebibliography}{99}
%-----------------------------------------------------------------------------------
%-----------------------------------------------------------------------------------


%-----------------------------------------------------------------------------------
%-----------------------------------------------------------------------------------


\bibitem{AlvesNovo}
C.O. Alves,
{\it \, Existence of positive solution for a nonlinear elliptic equation with saddle-like potential and nonlinearity with exponential critical growth in $\mathbb{R}^2$} To appear in Milan J. Math.

%-----------------------------------------------------------------------------------
%-----------------------------------------------------------------------------------



\bibitem{AdoOM} 
C. O. Alves, J. M. B. do \'O and O. H. Miyagaki, 
{\it On nonlinear perturbations of a periodic elliptic problem in $\mathbb{R}^2$ involving critical growth, } 
Nonlinear Anal. 56 (2004), 781--791.

%-----------------------------------------------------------------------------------
%-----------------------------------------------------------------------------------



%\bibitem{AOS}C.O. Alves, J.M. B. do \'{O} and M.A.S. Souto,
%{\it Local mountain-pass  for a class of elliptic problems involving critical growth.} 
%Nonlinear Anal. {\bf{46}} (2001),495-510.

%-----------------------------------------------------------------------------------
%-----------------------------------------------------------------------------------


\bibitem{AF2} 
C.O. Alves and G.M. Figueiredo, Multiplicity of Positive Solutions for a class of quasilinear problems in $\mathbb{R}^N$ via penalization method. Adv. Nonlinear Stud. {\bf 5},  (2005) 551-572.

%-----------------------------------------------------------------------------------
%-----------------------------------------------------------------------------------


%\bibitem{saddle} C.O. Alves and O. H. Miyagaki, {\it  A critical nonlinear fractional elliptic equation with  saddle-like potentical in $\mathbb{R}^N$,}  
%arXiv:1506.06533 

%-----------------------------------------------------------------------------------
%-----------------------------------------------------------------------------------


\bibitem{penalizacao} 
C.O. Alves and O. H. Miyagaki, 
{\it  Existence and concentration of solution for a class of fractional  elliptic equation in $\mathbb{R}^N$ via penalization method ,}  Preprint 

%-----------------------------------------------------------------------------------
%-----------------------------------------------------------------------------------



\bibitem{ARMPT} 
A. Ambrosetti and  P.H. Rabinowitz, 
{\it Dual variational methods in critical point theory and applications}, 
J. of Funct. Anal. {\bf 14}, (1973), 349-381.

%-----------------------------------------------------------------------------------
%-----------------------------------------------------------------------------------


\bibitem{barrios} 
B.\ Barrios, E.\ Colorado, A.\ de Pablo, U.\ S\'anchez, 
\textit{On some critical problems for the fractional Laplacian operator}, 
J. Differential  Equations {\bf 252} (2012), 613--6162.


%-----------------------------------------------------------------------------------
%-----------------------------------------------------------------------------------

\bibitem{BGe} R.M. \ Blumenthal and R.K.\ Getoor, Some theorems on stable processes, Trans. Amer. Math. Soc. 95 (1960), 263-273.

%-----------------------------------------------------------------------------------
%-----------------------------------------------------------------------------------



\bibitem{colorado} 
C. Br\"andle, E. Colorado and  U. S\'anchez,
{\it \,  A concave-convex elliptic problem involving the fractional Laplacian,} 
Proc. R. Soc. Edinb. A. {\bf 143} (2013), 39--71.


\bibitem{cabre} 
X. Cabr\'e and Y. Sire, 
{\it \, Nonlinear equations for fractional laplacians, I: Regularity,maximum principles, and Hamiltonian estimates,} 
Ann. Inst. H. Poincar\'e Non Lin\'eare {\bf 31} (2014), 23-53.

%-----------------------------------------------------------------------------------
%-----------------------------------------------------------------------------------

\bibitem{caffarelli} 
L.~Caffarelli and L.~Silvestre, 
\textit{An extension problems related to the fractional Laplacian,} 
Comm. Partial Differential Equation  {\bf 32} (2007), 1245--1260.

%-----------------------------------------------------------------------------------
%-----------------------------------------------------------------------------------

\bibitem{Cao}  D.M. Cao, {\it \,  Nontrivial solution of semilinear elliptic equation with critical exponent in $\R^2,$} Comm. Partial Differential Equation {\bf 17} (1992), 407-435.

\bibitem{CW} 
X. Chang and Z.Q.Wang, 
{\it \,  Nodal and multiple solutions of nonlinear problems involving the fractional Laplacian,} 
J. Differential Equations {\bf 256}(2014), 2965-2992.

%-----------------------------------------------------------------------------------
%-----------------------------------------------------------------------------------


\bibitem{Chen} 
G. Chen and Y. Zheng, 
{\it Concentration phenomena for fractional noninear Schr\"odinger equations,} 
Comm. Pure Appl. Anal. {\bf 13} (2014),  2359--2376.


%-----------------------------------------------------------------------------------
%-----------------------------------------------------------------------------------

\bibitem{Davila} 
J. D\'avila, M. del Pino and J. C. Wei, 
{\it  Concentrating standing waves for the fractional nonlinear Schr\"odinger equation,} 
J. Differential Equations {\bf 256} (2014), 858--892.

%-----------------------------------------------------------------------------------
%-----------------------------------------------------------------------------------


\bibitem{DMR} 
D.G.\ de Figueiredo, O.H.\ Miyagaki and  B.\ Ruf,
\textit{Elliptic equations in $\R^2$ with nonlinearities in the critical growth range},
Calc. Var. Partial Differential Equations {\bf 3} (1995), 139--153. 

%-----------------------------------------------------------------------------------
%-----------------------------------------------------------------------------------


\bibitem{Pino} 
M. del Pino and P.L. Felmer, 
{\it Local Mountain Pass for semilinear elliptic problems in unbounded domains}. Calc. Var. Partial Differential Equations  {\bf{4}} (1996), 121-137.

%-----------------------------------------------------------------------------------
%-----------------------------------------------------------------------------------



\bibitem{nezza} 
E.\ Di Nezza, G.\ Palatucci and E.\ Valdinoci, 
{\it Hitchhiker's guide to the fractional Sobolev spaces,}
Bull. Sci. Math. {\bf 136} (2012), 521--573.


%-----------------------------------------------------------------------------------
%-----------------------------------------------------------------------------------




\bibitem{DMS} 
J. M. do \'O, O. H. Miyagaki and  M. Squassina, 
\textit{Nonautonomous fractional problems  with  exponential growth}, to appear in NoDEA

%-----------------------------------------------------------------------------------
%-----------------------------------------------------------------------------------


\bibitem{OS} 
J. M. B.  do \'O and M.A.S. Souto, 
{\it On a class of nonlinear Schodinger equations in $\mathbb{R}^2$ involving critical growth}, 
J. Differential Equations  {\bf{174}} (2001), 289-311.


%-----------------------------------------------------------------------------------
%-----------------------------------------------------------------------------------



\bibitem{Moustapha} 
M. M. Fall, F. Mahmoudi and E. Valdinoci, 
{\it Ground states and concentration phenomena for the fractional Schr\"odinger equation,} 
arXiv:1411.0576v2[math.AP] 24 Apr 2015

%-----------------------------------------------------------------------------------
%-----------------------------------------------------------------------------------


\bibitem{FQT} 
P. Felmer, A Quass and J. Tan, 
{\it Positive solutions of nonlinear Schr\"odinger equation with the fractional Laplacian}, 
Proc. Roy. Soc. Edinburgh A {\bf 142} (2012), 1237–-1262.

%-----------------------------------------------------------------------------------
%-----------------------------------------------------------------------------------


%\bibitem{FW} A. Floer  and A. Weinstein, 
%{\it Nonspreading wave packets for the cubic Schr\"{o}dinger equations with bounded potential,} 
%J. Funct. Anal.  {\bf{69}} (1986), 397-408.


%-----------------------------------------------------------------------------------
%-----------------------------------------------------------------------------------


\bibitem{Frank} 
R. Frank and E. Lenzmann,
\textit{ Uniqueness of non-linear ground states for fractional laplacians in $\R,$} 
Acta Math. \textbf{210} (2013), 261--318.

%-----------------------------------------------------------------------------------
%-----------------------------------------------------------------------------------


%\bibitem{Frank1} R.L. Frank, E. Lenzmann and L. Silvestre, {\it  Uniqueness of radial solutions for the fractional laplacian,} %arXiv:1302.2652v2[math.AP]23 Mar 2015.

%-----------------------------------------------------------------------------------
%-----------------------------------------------------------------------------------


\bibitem{Antonio} 
A.\ Iannizzotto and  M.\ Squassina,  
$1/2$\textit{-laplacian problems with exponential nonlinearity,} 
J. Math. Anal. Appl. {\bf 414} (2014), 372--385.

%-----------------------------------------------------------------------------------
%-----------------------------------------------------------------------------------


%\bibitem{moser} 
%J.\ Moser, 
%\textit{A sharp form of an inequality by N.\ Trudinger}, 
%Indiana Univ. Math. J. {\bf 20} (1970), 1077--1092. 

%-----------------------------------------------------------------------------------
%-----------------------------------------------------------------------------------


%\bibitem{Oh1} 
%Y.J. Oh, 
%{\it Existence of semi-classical bound states of nonlinear Schr\"{o}dinger equations with potentials on the class 	$(V)_a$,}  
% Comm. Partial Differential Equations  {\bf{13}} (1988), 1499-1519.

%-----------------------------------------------------------------------------------
%-----------------------------------------------------------------------------------


\bibitem{Ozawa} 
T.~Ozawa, 
\textit{On critical cases of Sobolev's inequalities,}
J. Funct. Anal. \textbf{127} (1995), 259--269. 

%-----------------------------------------------------------------------------------
%-----------------------------------------------------------------------------------



%\bibitem{R} 
%P.H. Rabinowitz,  
%{\it On a class of nonlinear Schr\"{o}dinger equations}, 
%Z. Angew Math. Phys. {\bf{43}} (1992),270-291.

%-----------------------------------------------------------------------------------
%-----------------------------------------------------------------------------------


\bibitem{Secchi} 
S. Secchi, 
\textit{Ground state solutions for nonlinear fractional Schr\"odinger equations in $\Re^N$,}  
J. Math. Phys.   {\bf 54}  (2013), 031501-17 pages.

%-----------------------------------------------------------------------------------
%-----------------------------------------------------------------------------------


\bibitem{ShangZhang} 
X. Shang and J. Zhang, 
{\it Concentrating solutions of nonlinear fractional Schr\"odinger equation with potentials,} 
J. Differential Equations {\bf 258} (2015), 1106--1128.

%-----------------------------------------------------------------------------------
%-----------------------------------------------------------------------------------


\bibitem{Shang} 
X. Shang, J. Zhang and  Y. Yang, 
\textit{On fractional Sch\"odinger equation in $\Re^N$ with critical growth,} 
J. Math. Phys. {\bf 54} (2013), 121502-19 pages.

%-----------------------------------------------------------------------------------
%-----------------------------------------------------------------------------------


%\bibitem{W}X. Wang, 
%{\it On concentration of positive bound states of nonlinear Schr$\ddot{o}$dinger equations,} 
%Comm. Math. Physical {\bf{53}} (1993), 229-244.

%-----------------------------------------------------------------------------------
%-----------------------------------------------------------------------------------


\bibitem{willem}  
M. Willem,\textit{Minimax Theorems,} 
Progress in Nonlinear Differential Equations and their Applications, 24, Birkh$\ddot{a}$user 1996. 

%-----------------------------------------------------------------------------------
%-----------------------------------------------------------------------------------


%-----------------------------------------------------------------------------------
%-----------------------------------------------------------------------------------
\end{thebibliography}
\end{document}